\documentclass[11pt]{amsart}

\usepackage{mathrsfs, amsmath, amsthm, amsfonts, amssymb, amscd, graphicx, epstopdf, color,hyperref}
\usepackage[all]{xy}
\usepackage{pinlabel}

\newtheorem*{thm*}{Main Theorem}
\newtheorem*{conj*}{Generalized Property $\bs{R}$ Conjecture}

\theoremstyle{definition}
\newtheorem*{thanks*}{Acknowledgements}
        
\def\Z{\mathbb{Z}}
\def\S{\mathbb{S}}

\def\N{\mathbb{N}}
\def\bs{\boldsymbol}

\def\a{\alpha}

\def\g{\gamma}

\def\d{\partial}

\def\<{\langle}
\def\>{\rangle}

\def\cl{\overline}

\begin{document}

\title{Handle Number One Links and Generalized Property $R$}
\author{Michael J. Williams}
\address{\hskip-\parindent
	Department of Mathematics\\
	University of California\\
	Santa Barbara, CA 93106\\
	USA}
\email{mikew@math.ucsb.edu}
\thanks{Research supported by a UCOP Postdoctoral Fellowship.}

\begin{abstract} 
It is shown that if the exterior of a link $L$ in the three sphere admits a genus 2 Heegaard splitting, then $L$ has Generalized Property $R$. 
\end{abstract}	

\maketitle

\section{Introduction} 

A knot $K$ in the three sphere $\S^3$ is said to have \textbf{Property $\bs{R}$} if $K$ is the unknot, or $0$-framed Dehn surgery on $K$ does not yield $\S^1 \times \S^2$. Gabai \cite{gabai:foliations3} proved that every knot in $\S^3$ has Property $R$. A generalization to links is the following 

\begin{conj*}
If Dehn surgery on an $n$-component framed link $L \subset \S^3$ yields $\#_n(\S^1 \times \S^2)$, then there is a sequence of handleslides which converts $L$ to the $n$-component $0$-framed unlink. 
\end{conj*}

Recall that the \textit{handleslide} is one of the basic operations of the Kirby Calculus in which a framed link component gets replaced by its band connected sum with the framing of another framed component; see \cite[Chapter 4]{gs:4} for more details. All the framings on a link $L$ satisfying the hypothesis of Generalized Property $R$ must be the $0$-framing, and the components of $L$ must be algebraically unlinked; this is due to classical homological arguments (see \cite[Proposition 2.2]{st:fibered}). We will say that a link $L$ has \textbf{Generalized Property $\bs{R}$} if $L$ satisfies the the Generalized Property $R$ Conjecture. 

Part of the significance of the Generalized Property $R$ Conjecture lies in the realm of smooth $4$-dimensional handle structures. The closed 3-manifold $M$ obtained by integral Dehn surgery on an $n$-component link $L$ in $\S^3$ is the boundary of a compact, orientable, simply connected 4-manifold $X$ obtained by attaching 2-handles to the standard smooth 4-ball $B^4$ along $L$ with framings given by the Dehn surgery slopes; see \cite[Chapters 4 and 5]{gs:4} for details. If $M = \#_n(\S^1 \times \S^2)$, then the $\S^2$ factors of the $\S^1 \times \S^2$ summands serve as attaching regions for attaching $3$-handles to $X$; this yields a smooth homotopy $4$-ball $X'$ with $\d X' = \S^3$. If $L$ has Generalized Property $R$, then the $2$-handles and $3$-handles can be arranged to cancel in pairs, giving a standard handle structure on $B^4$. Hence $X'$ is diffeomorphic to $B^4$, preventing $X'$ from generating a counterexample to the Smooth Poincar\'e Conjecture in dimension $4$. For more information, see \cite[Section 2]{scharlemann:generalized} and \cite[Section 3]{fgmw:man}.

A \textbf{tunnel number one link in $\bs{\S^3}$} is a link whose exterior is not a handlebody, yet admits a genus 2 Heegaard splitting which decomposes the exterior into a handlebody and a compression body. Therefore, the exterior of such a link is obtained by attaching a $2$-handle to a genus 2 \textit{handlebody} $H$ along a simple closed curve on $\d H$.  Proposition 3.1 of \cite{st:fibered} asserts that the only $2$--component tunnel number one link to admit an integral $\#_2(\S^1 \times \S^2)$ surgery is the $2$-component unlink. Hence tunnel number one links satisfy the Generalized Property $R$ Conjecture. In this article, we address the conjecture for a larger class of links.  

A \textbf{handle number one link in $\bs{\S^3}$} is a link whose exterior is not a compression body, yet admits a genus 2 Heegaard splitting which decomposes the exterior into compression bodies; this terminology first appeared in \cite{sedgwick:genus}; also see \cite{kobayashi:scharlemann} and \cite{ms:closed} for examples. This generalizes tunnel number one links only insofar as we do not require the Heegaard surface to bound a handlebody on one side. It is easy to see that a handle number one link can have up to four components. It is straightforward to see that a link has handle number one if and only if its exterior has Heegaard genus $2$. Our main result below is inspired by \cite[Proof of Proposition 3.1]{st:fibered}, and resolves the Generalized Property $R$ Conjecture for all links with Heegaard genus $2$.

\begin{thm*}
Suppose that $L \subset \S^3$ is an $n$--component link $(n > 1)$ whose exterior admits a genus $2$ Heegaard splitting. If $0$--framed surgery on $L$ yields $\#_n(\S^1 \times \S^2)$, then 
\begin{itemize}
\item $n=2$, and
\item the components of $L$ are unknots.
\end{itemize}
Consequently $L$ has Generalized Property $R$.
\end{thm*}

Our notation for Dehn surgery will be the following. Given a link $L \subset \S^3$ with components $L_1, \dots, L_n$ and corresponding framings $r_1, \dots, r_n \in \Z$, let $L(r_1,\dots,r_n)$ denote the \textit{$3$--manifold obtained by Dehn surgery on the framed link $L$}. If all the $r_i$ are equal to a particular $r \in \Z$, we will refer to $L$ as an $r$--framed link. 

\begin{thanks*}
I would like to thank Martin Scharlemann for helpful conversations.
\end{thanks*}\label{sec:introduction}
\section{Proof of the Main Theorem}

The proof will involve fundamental group calculations. In order to keep our notation as simple as possible, we will adopt the following convention: if $M$ is a connected $3$-manifold with basepoint $p \in M$, and $\g \subset M$ is a loop, then we will also let $\g$ denote the homotopy class of $\g$ in $\pi_1(M,p)$, possibly under an appropriate change-of-basepoint isomorphism when $p \notin \g$. If $G$ is a group and $S \subset G$, let $\<\<S\>\>$ denote the \textbf{normal closure} of $S$ in $G$. 

\begin{proof}[Proof of the Main Theorem]
Let $L \subset \S^3$ be an $n$--component link $n > 1$ whose exterior $E(L)$ admits a genus $2$ Heegaard splitting. Therefore, any Dehn surgery on $L$ admits a genus $2$ Heegaard splitting. By additivity of Heegaard genus under connected sums, we see that $\#_m(\S^1 \times \S^2)$ has Heegaard genus $m$ for any $m \in \N$. By assumption, $0$--framed Dehn surgery on $L$ yields $\#_n(\S^1 \times \S^2)$; so $n=2$.

Let $L_1$ and $L_2$ be the components of $L$ with respective regular neighborhoods $N(L_1)$ and $N(L_2)$. So $L$ has exterior $E(L)=\cl{\S^3 - (N(L_1) \cup N(L_2))}$. Let $\d_1$ and $\d_2$ be the boundary components of $E(L)$, so that $\d_i$ corresponds to $L_i$. Let $m_i,l_i$ be a standard oriented meridian-longitude pair for $\d_i$ for $i=1,2$. 

Now, if $L$ has tunnel number one, then we are done by \cite[Proposition 3.1]{st:fibered}. So we may assume that there is a genus 2 Heegaard splitting of $E(L)$ which separates the boundary components. Hence the exterior of $L$ can be realized as $E(L)=W \cup_\a (\text{2-handle})$ where $W$ is a compression body with $\d_-W=\d_1$, genus($\d_+W$)=2, and $\a$ is a nonseparating curve on $\d_+W$. Note that $W$ is obtained by attaching a 1-handle to $(\text{torus}) \times [0,1]$ on $(\text{torus}) \times \{1\}$. If we slide the endpoints of the core of the 1-handle together, we obtain a simple closed curve $\tau$ which meets $(\text{torus}) \times \{1\}$ in one point; we push $\tau$ toward $\d_-W=(\text{torus}) \times \{0\}=\d_1$ to meet it in a point $p$. Give $\tau$ an arbitrary orientation and situate the curves $m_1$ and $l_1$ on $\d_1$ so that $m_1 \cap l_1 = p$; see Figure~\ref{fig:compressionbody1}. It is now clear that 
\begin{align*}
&\pi_1(W,p) \cong \< m_1,l_1,\tau : m_1l_1=l_1m_1 \> \cong (\Z \oplus \Z) * \Z \ , \\
&\pi_1(E(L),p) \cong  \pi_1(W,p)/\<\< \a \>\> \ , \ \text{and} \\
&\pi_1(L(0,0),p) \cong \pi_1(E(L),p)/\<\<l_1,l_2\>\> \cong \#_2(\S^1 \times \S^2) \cong \Z * \Z \ .
\end{align*}

\begin{figure}[htbp]
  \centering
   \labellist
   \small \hair 2pt
   \pinlabel {$\d_1$} at 520 425
   \pinlabel {$\color{blue}{\bs{l_1}}$} at 505 212
   \pinlabel {$\color{red}{\bs{m_1}}$} at 480 290
   \pinlabel {$p$} [r] at 338 293
   \pinlabel {$\color{green}{\bs{\tau}}$} at 230 440
   \pinlabel {$W$} at 290 077
   \endlabellist  
\includegraphics[scale=.32]{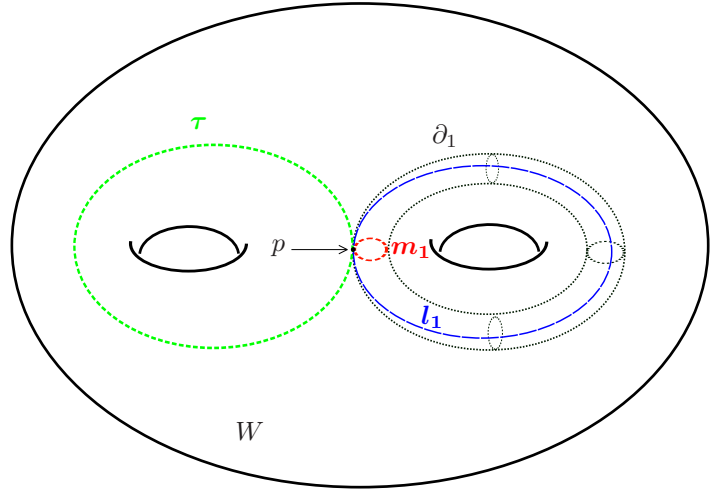}
\caption{The compression body $W$ with curves $m_1$ and $l_1$ on $\d_-W=\d_1$ and the ``tunnel curve" $\tau$. The curves $m_1$, $l_1$, and $\tau$ all meet the basepoint $p$. Note that $\tau \cup \d_1$ forms a spine for $W$.}
\label{fig:compressionbody1}
\end{figure}

To simplify notation, the basepoint $p \in W \subset E(L) \subset L(0,0)$ will be suppressed in fundamental group expressions for the remainder of the proof. We have epimorphisms
$$\xymatrix{\pi_1(W) \ar@{>>}[r]^{\phi \ \ } & \pi_1(E(L)) \ar@{>>}[r]^{\rho \ \ } & \pi_1(L(0,0))},$$
where $\phi$ is reduction mod $\<\< \a \>\>$, and $\rho$ is reduction mod $\<\< l_1,l_2 \>\>$. Consider another sequence of epimorphisms
$$\xymatrix{\pi_1(W) \ar@{>>}[r]^{\phi \ \ } & \pi_1(E(L)) \ar@{>>}[r]^{\rho_0 \ \ \ } & \pi_1(E(L))/\<\< l_1 \>\>},$$
where $\rho_0$ is reduction mod $\<\< l_1 \>\>$. Note that the epimorphism $\rho$ \textit{factors through $\rho_0$}, that is, there is an epimorphism 
$$\xymatrix{\pi_1(E(L))/\<\<l_1\>\> \ar@{>>}[r]^{\ \ \ \rho_1} & \pi_1(L(0,0))},$$
namely reduction mod $\<\< l_2 \>\>$, so that $\rho=\rho_1 \circ \rho_0$. Since $(\rho_0 \circ \phi) (l_1) = 1$, the epimorphism $\rho \circ \phi$ naturally descends to an epimorphism
$$\xymatrix{\pi_1(W)/\<\< l_1 \>\> \ar@{>>}[r]^{\phi' \ \ } & \pi_1(E(L))/\<\< l_1 \>\>}.$$
Note that $W \cup_{l_1} (\text{2-handle})$ is a punctured genus $2$ handlebody; so $\pi_1(W)/\<\< l_1 \>\> \cong \Z * \Z$. 
In summary, we have a commutative diagram of epimorphisms
$$\xymatrix{
(\Z \oplus \Z)*\Z \ar@{=}[r] & \pi_1(W) \ar@{>>}[r]^{\phi \ \ \ \ } \ar@{>>}[d] & \pi_1(E(L)) \ar@{>>}[d]_{\rho_0} \ar@{>>}[r]^{\rho \ \ } & \pi_1(L(0,0)) \ar@{=}[d] \\
\Z*\Z \ar@{=}[r] & \pi_1(W)/ \<\<l_1\>\> \ar@{>>}[r]^{\phi' \ \ }  & \pi_1(E(L))/\<\<l_1\>\> \ \ar@/__/@{>>}[ur]_{\ \ \rho_1} & \Z*\Z}.$$
Since the group $\Z * \Z$ is Hopfian (\textit{every self-epimorphism is an isomorphism}), we must have all isomorphisms in the sequence
$$\xymatrix{
\Z*\Z \ar@{=}[r] & \pi_1(W)/\<\<l_1\>\> \ar@{>>}[r]^{\phi' \ \ } & \pi_1(E(L))/\<\<l_1\>\> \ar@{>>}[r]^{\ \ \rho_1} & \pi_1(L(0,0)) \ar@{=}[r] & \Z*\Z}.$$
We now see that $\pi_1(E_1) \cong \pi_1(E(L))/\<\<l_1\>\> \cong \Z * \Z$, where $E_1$ is the $3$--manifold obtained from $E(L)$ Dehn filling along $l_1 \subset \d N(L_1)$. By \cite[Theorem 5.2]{hempel:3}, the Prime Decomposition Theorem, and the fact that genus $2$ $3$--manifolds satisfy the Poincar\'e Conjecture (see \cite{bm:smith}), we conclude that $$E_1 \cong (\S^1 \times D^2) \# (\S^1 \times \S^2) \ .$$ Furthermore, Dehn filling on $E_1$ along $l_2 \subset \d N(L_2)$ gives $L(0,0) \cong (\S^1 \times \S^2) \# (\S^1 \times \S^2)$, so $l_2$ corresponds to a meridian curve in the $\S^1 \times D^2$ connected-summand of $E_1$. Note that this is just the topological realization of the the epimorphisms $\rho_0$ and $\rho_1$. It follows that $m_2$ generates a free factor of $\pi_1(E_1)$. Hence $\pi_1(E_1)/\<\<m_2\>\> \cong \Z$. 

Let $L_1(0)$ denote the $3$--manifold obtained by $0$--framed Dehn surgery on the knot $L_1$. We see that   
\begin{align*}
\pi_1(L_1(0)) &\cong \pi_1(E(L))/\<\< l_1,m_2 \>\> \\
&\cong \bigl( \pi_1(E(L))/\<\< l_1 \>\> \bigr) / \<\< m_2 \>\> \\
&\cong \pi_1(E_1)/\<\<m_2\>\> \\
&\cong \Z \ .
\end{align*}
By \cite[Theorem 5.2]{hempel:3}, we have that $L_1(0) \cong \S^1 \times \S^2$. Applying \cite[Corollary 8.3]{gabai:foliations3} establishes that $L_1$ must be the unknot. We can similarly establish that $L_2$ is also the unknot; this is accomplished by simply interchanging the roles of $L_1$ and $L_2$ throughout the proof. Finally, the result \cite[Proposition 3.2]{st:fibered} asserts that any $2$--component link containing an unknot has Generalized Property $R$. This completes the proof. 
\end{proof}
\label{sec:main_thm}

\bibliographystyle{amsplain}
\bibliography{H1GPR}

\end{document}